\documentclass[12pt]{amsart}
\usepackage{amssymb}

\usepackage[dvips]{graphicx}
\usepackage[dvips]{graphics}
\newtheorem{theorem}{{\bf Theorem}}
\newtheorem{prop}{{\bf Proposition}}
\newtheorem{cor}{Corollary}%
\newtheorem{rem}{Remark}

\begin{document}

\title{REPRESENTATIONS  OF BRAID GROUPS AND GENERALISATIONS}

\author[Bardakov]{Valerij G. Bardakov}
\address{Sobolev Institute of Mathematics, Novosibirsk 630090, Russia}
\email{bardakov@math.nsc.ru}

\author[Bellingeri]{Paolo Bellingeri}
\address{
Universit\`a  Milano-Bicocca, Dipartimento Matematica e Applicazioni, 20126 Milano, Italy}
\email{paolo.bellingeri@unimib.it}

\date{\today}

\begin{abstract}
We define and study  extensions of Artin's representation and braid monodromy representation 
to the case of topological and algebraical generalisations of braid groups. 
In particular we provide  faithful representations  of braid groups of oriented surfaces 
with boundary components as (outer) automorphisms of free groups.
We give also similar representations for  braid groups of non oriented surfaces 
with boundary components and
 we show a representation of braid groups of closed surfaces as outer automorphisms of free groups.
Finally, we provide faithful representations  of Artin-Tits groups  of type $\mathcal{D}$
as automorphisms of free groups.
\end{abstract}

\maketitle
\section{Introduction}
Let $F_n$ be the free group of rank $n$ with the set of generators $\{  x_1, x_2, \ldots, x_n \}$.
Assume further that ${\rm Aut}(F_n)$ is the automorphism group of $F_n$.
The Artin braid group $B_n$ can be represented as a subgroup of ${\rm Aut}(F_n)$.
This representation, due to   Artin himself, is defined associating to
any generator $\sigma_i$, for $i=1,2,\ldots,n-1$, of $B_n$ the following
automorphism of $F_n$:
$$
\sigma_{i} : \left\{
\begin{array}{ll}
x_{i} \longmapsto x_{i} \, x_{i+1} \, x_i^{-1}, &  \\ x_{i+1} \longmapsto
x_{i}, & \\ x_{l} \longmapsto x_{l}, &  l\neq i,i+1.
\end{array} \right.
$$

Moreover (see for instance~\cite[Theorem 5.1]{LH}),
 any automorphism
$\beta $ of ${\rm Aut}(F_n)$ corresponds to an element of  $B_n$
if and only if $\beta $ satisfies  the following conditions:
\begin{eqnarray*}
& i) & \beta(x_i) = a_i^{-1} \, x_{s(i)} \, a_i,~~1\leq i\leq n,\\
& ii)& \beta(x_1x_2 \ldots x_n)=x_1x_2 \ldots x_n,
\end{eqnarray*}

where $s$ is a permutation from the symmetric group $S_n$ and $a_i\in F_n$.

Generalisations of Artin's representation have been provided  by Wada~\cite{W}
and further by Crisp and Paris~\cite{CP0}, in order to construct group invariants 
of oriented links.

Another interesting representation of $B_n$ is the \emph{braid monodromy representation} of 
$B_n$ into ${\rm Aut}(F_{n-1})$ (see last Section), that was proven to be faithful in~\cite{CP} and~\cite{PV}.

In this paper, we extend Artin's and braid monodromy representations to some generalisations of braid groups.

In the case of braid groups of oriented
surfaces with boundary components our representations
are faithful  (Theorem~\ref{theorem3annexe}). Moreover, in the case of surfaces of genus $g \ge 1$,
the induced representations
of  surface braid groups  as  outer automorphisms  hold faithful (Theorem~\ref{out}). 

In the case of closed surfaces,
as earlier remarked in~\cite{Bir} for braid groups of the sphere,
we cannot extend Artin's representation and we provide a representation in  the outer automorphism group
of a finitely generated free group. 

In the last Section we consider the braid monodromy representation of $B_n$ in ${\rm Aut}(F_{n-1})$ and  we provide a faithful representation of the $n$-th Artin-Tits group of type $\mathcal{D}$ in ${\rm Aut}(F_{n})$ (Proposition~\ref{thmartind}). 
\vspace{5pt}

\noindent{\bf Acknoledgments.} 
The research of the first author  has been supported by the Scientific
and Research Council of Turkey (T\"UB\.ITAK). The first author would like
to thank Prof. Hasan G\"umral, Prof. Vladimir Tolstykh
and other members of the Department of Mathematics
of Yeditepe University (Istanbul) for their
kind hospitality. We thank Prof.~Warren Dicks for pointing out  a doubtful 
argument on a previous version of Proposition~\ref{thmartind} and for 
suggesting us a different  approach for  an easier proof.

\section{Braid groups of orientable surfaces with boundary components} \label{firstsection}

\noindent \textbf{Surface braids as collections of paths.} Let $\mathcal{P}=\{p_1, \dots,
p_n\}$ be a set of $n$ distinct points (\emph{punctures}) in the
interior of a surface $\Sigma$.

 A \emph{geometric braid} on $\Sigma$ based at $\mathcal{P}$
is a collection $(\psi_1, \dots, \psi_n)$ of $n$ disjoint paths
(called \emph{strands}) on $\Sigma \times [0, 1]$ which run
monotonically with $t \in [0, 1]$ and such that $\psi_i(0)=(p_i,
0)$ and $\psi_i(1) \in \mathcal{P} \times \{1 \}$. Two braids are
considered to be equivalent if they are isotopic relatively to
the base points. The usual product of paths defines a group
structure on the equivalence classes of braids. This group,
denoted usually by $B_n(\Sigma)$, does not depend on the choice
of $\mathcal{P}$ and it is called  braid group on $n$ strands of
$\Sigma$. The $n$th braid group of the disk $D^2$, $B_n(D^2)$, is isomorphic to $B_n$.

In the following we will denote by $B_n(\Sigma_{g,p})$ the braid group on $n$ strands of
an orientable surface of genus $g$
with $p$  boundary components (we set $\Sigma_g =\Sigma_{g,0}$) and by
$B_n(N_{g,p})$ the braid group on $n$ strands of a non-orientable surface of genus $g$
with $p$ boundary components (we set $N_g = N_{g,0}$).

In this section we will consider  an orientable surface $\Sigma_{g,p}$ of genus $g \geq 0$
and with $p>0$ boundary
components. We set also $n\ge 2$.

 We denote by $\sigma_1,...,\sigma_{n-1}$ the standard
generators of the braid group $B_n$. Since $p>0$, we can embed a disk in
$\Sigma_{g,p}$  and therefore we can consider  $\sigma_1,...,\sigma_{n-1}$
 as elements of $B_n(\Sigma_{g,p})$.  Let also
$a_1,...,a_{g}, b_1,...,b_g, z_1,...,z_{p-1}$ be the generators
of $\pi_1(\Sigma_{g,p})$, where $z_i$'s denote loops around the holes.
Assume that the base point of the fundamental group is
the startpoint of the first strand. Then each element $\gamma\in
\pi_1(\Sigma_{g,p})$ determines an element denoted also by $\gamma$ in
$B_n(\Sigma_{g,p})$, by considering the braid whose first strand is describing the
curve $\gamma$ and other strands are constant.

Let us set $x_{2k-1}=a_k$ and $x_{2k}=b_k$  for $k=1, \dots, g$. 
According to  \cite{Bel} we have
that the group $B_n(\Sigma_{g, p})$  admits a presentation with  generators:
$$
\sigma_1, \ldots , \sigma_{n-1}, x_1, \ldots , x_{2g}, z_1, \ldots, z_{p-1},
$$
and  defining relations:

\noindent \noindent -- Braid relations:

$\sigma_i \sigma_{i+1} \sigma_i = \sigma_{i+1} \sigma_i \sigma_{i+1},~~1 \leq i \leq n-2,$

$\sigma_i \sigma_{j} = \sigma_{j} \sigma_i,~~|i - j| > 1,~~1 \leq i,j \leq n-1,$

\noindent -- Mixed relations:

\noindent (R1)~$x_r \sigma_i = \sigma_i x_r,$~~$i \neq 1,$~~$1 \leq r \leq 2g,$

\noindent (R2)~$(\sigma_1^{-1} x_r \sigma_1^{-1}) x_r = x_r (\sigma_1^{-1} x_r \sigma_1^{-1}),$
~~$1 \leq r \leq 2g,$

\noindent (R3)~$(\sigma_1^{-1} x_s \sigma_1) x_r = x_r (\sigma_1^{-1} x_s \sigma_1)$, $1 \leq s < r \leq 2g$, $(s, r) \neq (2m-1, 2m)$,

\noindent (R4)~$(\sigma_1^{-1} x_{2m-1} \sigma_1^{-1}) x_{2m} = x_{2m} (\sigma_1^{-1} x_{2m-1}
\sigma_1),$
~~$1 \leq m \leq g,$

\noindent (R5)~$z_j \sigma_i = \sigma_i z_j,$~~$i \neq 1,$~~$1 \leq j \leq p-1,$

\noindent (R6)~$(\sigma_1^{-1} z_j \sigma_1) x_r = x_r (\sigma_1^{-1} z_j \sigma_1),$
~~$1 \leq  r \leq 2g,$~~$1 \leq j \neq p-1,$

\noindent (R7)~$(\sigma_1^{-1} z_j \sigma_1) z_l = z_l (\sigma_1^{-1} z_j \sigma_1),$
~~$1 \leq  j < l \leq p-1,$

\noindent (R8)~$(\sigma_1^{-1} z_j \sigma_1^{-1}) z_j = z_j (\sigma_1^{-1} z_j \sigma_1),$
~~$1 \leq  j \leq p-1.$

\medskip

Associating to any surface braid the corresponding permutation
one obtains  a surjective homomorphism
$$
\pi : B_n(\Sigma_{g,p}) \longrightarrow S_n,
$$
such that
$\pi(\sigma_i) = (i, i+1), ~i = 1, \ldots, n-1$,
$\pi(x_r) = \pi (z_j)=e$ for $1 \leq  r \leq 2g$ and $1 \leq  j \leq p-1
$.

In \cite{Bel} the second author considered the subgroup $D_n(\Sigma_{g,p}) = \pi^{-1}(S_{n-1})$ and found its
generators. We provide a set of defining relations of $D_n(\Sigma_{g,p})$
using the well-known  Reidemeister-Schreier's method (see \cite[Chap. 2]{KMS}).

Let  $M_n = \{ m_l ~|~1 \leq l \leq n \}$ be the set  defined as follows:
$$
m_l = \sigma_{n-1} \ldots \sigma_l,~l = 1, \ldots, n-1,~m_n=1.
$$
It is easy to prove (see~\cite{Bel}) that   $|B_n(\Sigma_{g,p}) : D_n(\Sigma_{g,p})| = n$ and that
 $M_n$ is a Schreier set of coset representatives of
$D_n(\Sigma_{g,p})$ in $B_n(\Sigma_{g,p}).$
Define the map $^- : B_n(\Sigma_{g,p}) \longrightarrow M_n$ which takes an element
$w \in B_n(\Sigma_{g,p})$
into the representative $\overline{w}$ from $M_n$. The element
$w \overline{w}^{-1}$ belongs to $D_n(\Sigma_{g,p})$ and,  by Theorem 2.7 from  \cite{KMS},
the group $D_n(\Sigma_{g,p})$ is generated by
$$
s_{\lambda, a} = \lambda a \cdot (\overline{\lambda a})^{-1},
$$
where $\lambda$ runs over the set $M_n$ and $a$ runs over the set of generators of
$B_n(\Sigma_{g,p})$.
\\

\noindent {\bf Case 1.} If $a \in \{ \sigma_1, \ldots, \sigma_{n-1} \}$, then we find the generators
$$
\tau_k = \sigma_{n-1} \ldots \sigma_{k+1} \sigma_k^2 \sigma_{k+1}^{-1} \ldots \sigma_{n-1}^{-1},~
k=1,\ldots,n-2,~\tau_{n-1}=\sigma_{n-1}^2.
$$
\\

\noindent {\bf Case 2.} If $a \in \{ x_1, \ldots, x_{2g} \}$, then we find the generators
$$
w_r = \sigma_{n-1} \ldots \sigma_{1} x_r \sigma_{1}^{-1} \ldots \sigma_{n-1}^{-1},~
r=1,\ldots,2g.
$$
\\

\noindent {\bf Case 3.} If $a \in \{ z_1, \ldots, z_{p-1} \}$, then we find the generators
$$
\xi_j = \sigma_{n-1} \ldots \sigma_{1} z_j \sigma_{1}^{-1} \ldots \sigma_{n-1}^{-1},~
j=1,\ldots,p-1.
$$
\\

To find  defining relations of $D_n(\Sigma_{g,p})$ we define
a rewriting process $\tau $. It allows us to rewrite a word which is written in the generators
of $B_n(\Sigma_{g,p})$ and to present an element in $D_n(\Sigma_{g,p})$ as a word in the generators
of $D_n(\Sigma_{g,p})$.
Let us associate to the reduced word
$$
u = a_1^{\varepsilon_1} \, a_2^{\varepsilon_2} \ldots
a_{\nu}^{\varepsilon_{\nu}},~~\varepsilon_l = \pm 1,
$$
where
$$
a_l \in
\{\sigma_1, \sigma_2, \ldots, \sigma_{n-1},~~
 x_1, x_2, \ldots, x_{2g},~~
z_1, z_2, \ldots, z_{p-1}
\},
$$
the word
$$
\tau(u) = s_{k_1,a_1}^{\varepsilon_1} \,  s_{k_2,a_2}^{\varepsilon_2}
\ldots s_{k_{\nu},a_{\nu}}^{\varepsilon_{\nu}}
$$
in the generators of $D_n(\Sigma_{g,p})$, where $k_j$ is a representative of the ($j-1$)th
initial segment
of the word $u$ if $\varepsilon_j = 1$ and $k_j$ is a representative of the $j$th
initial segment of
the word $u$ if
$\varepsilon_j = -1$.
By \cite[Theorem 2.9]{KMS}, the group $D_n(\Sigma_{g,p})$ is defined by relations
$$
r_{\mu,\lambda} = \tau (\lambda  \, r_{\mu} \,  \lambda^{-1}),~~\lambda \in
M_n,
$$
where $r_{\mu}$ is a defining relation of $B_n(\Sigma_{g,p})$.

\begin{prop} \label{propdn1}
The group $D_n(\Sigma_{g,p})$  admits a presentation with the generators
\\

$
\sigma_1, \ldots, \sigma_{n-2}, x_1, \ldots, x_{2g}, z_1, \ldots, z_{p-1},
\tau_1, \ldots, \tau_{n-1},~w_1, \ldots, w_{2g},~\xi_1, \ldots,$
$
\xi_{p-1};
$

and   relations:
\\

\noindent -- Braid relations

\noindent (B1)~$\sigma_i \sigma_{i+1} \sigma_i = \sigma_{i+1} \sigma_i \sigma_{i+1},~~1 \leq i \leq n-3,$

\noindent (B2)~$\sigma_i \sigma_{j} = \sigma_{j} \sigma_i,~~|i - j| > 1,~~1 \leq i,j \leq n-2,$

\noindent (B3)~$\sigma_k^{-1} \tau_l \sigma_k = \tau_l,~k \neq l-1, l,$

\noindent (B4)~$\sigma_{l-1}^{-1} \tau_l \sigma_{l-1} = \tau_{l-1},$

\noindent (B5)~$\sigma_l^{-1} \tau_l \sigma_l = \tau_l \tau_{l+1} \tau_l^{-1},~l \neq n-1.$

\noindent -- Mixed relations

\noindent (R1.1)~$x_r \sigma_i = \sigma_i x_r,$~~$2 \leq i \leq n-2,$~~$1 \leq r \leq 2g,$

\noindent (R1.2)~$x_r \tau_i = \tau_i x_r,$~~$2 \leq i \leq n-1,$

\noindent (R1.3)~$w_r \sigma_i = \sigma_i w_r,$~~$1 \leq r \leq 2g,$~~$ \leq i \leq n-2,$

\noindent (R2.1)~$(\sigma_1^{-1} x_r \sigma_1^{-1}) x_r = x_r (\sigma_1^{-1} x_r \sigma_1^{-1}),$
~~$1 \leq r \leq 2g,$

\noindent (R2.2)~$x_r^{-1} w_r x_r = \tau_1^{-1} w_r \tau_1,$
~~$1 \leq r \leq 2g,$

\noindent (R2.3)~$x_r^{-1} \tau_1 x_r = \tau_1^{-1} w_r \tau_1 w_r^{-1} \tau_1,$
~~$1 \leq r \leq 2g,$

\noindent (R3.1)~$(\sigma_1^{-1} x_s \sigma_1) x_r = x_r (\sigma_1^{-1} x_s \sigma_1)$, $1 \leq s < r \leq 2g$, $(s, r) \neq (2m-1, 2m),$

\noindent (R3.2)~$x_r^{-1} (\tau_1^{-1} w_s \tau_1) x_r = \tau_1^{-1} w_s \tau_1$, $1 \leq s < r \leq 2g$, $(s, r) \neq (2m-1, 2m)$,

\noindent (R3.3)~$x_s w_r = w_r x_s,$
~~$1 \leq s < r \leq 2g,$~~$(s, r) \neq (2m-1, 2m),$

\noindent (R4.1)~$(\sigma_1^{-1} x_{2m-1} \sigma_1^{-1}) x_{2m} = x_{2m} (\sigma_1^{-1} x_{2m-1}
\sigma_1),$~~$1 \leq m \leq g,$

\noindent (R4.2)~$x_{2m}^{-1} ( \tau_1^{-1} w_{2m-1} ) x_{2m} = \tau_1^{-1} w_{2m-1} \tau_1,$
~~$1 \leq m \leq g,$

\noindent (R4.3)~$x_{2m-1}^{-1} w_{2m}  x_{2m-1} = \tau_1^{-1} w_{2m},$
~~$1 \leq m \leq g,$

\noindent (R5.1)~$z_j \sigma_i = \sigma_i z_j,$~~$2 \leq i \leq n-2,$~~$1 \leq j \leq p-1,$

\noindent (R5.2)~$z_j \tau_i = \tau_i z_j,$~~$2 \leq i \leq n-1,$~~$1 \leq j \leq p-1,$

\noindent (R5.3)~$\xi_j \sigma_i = \sigma_i \xi_j,$~~$1 \leq j \leq p-1,$~~$1 \leq i \leq n-2,$

\noindent (R6.1)~$(\sigma_1^{-1} z_j \sigma_1) x_r = x_r (\sigma_1^{-1} z_j \sigma_1),$
~~$1 \leq  r \leq 2g,$~~$1 \leq j \leq p-1,$

\noindent (R6.2)~$x_r^{-1} (\tau_1^{-1} \xi_j \tau_1) x_r = \tau_1^{-1} \xi_j \tau_1,$
~~$1 \leq  r \leq 2g,$~~$1 \leq j \leq p-1,$

\noindent (R6.3)~$z_j w_r = w_r z_j,$
~~$1 \leq  r \leq 2g,$~~$1 \leq j \leq p-1,$

\noindent (R7.1)~$(\sigma_1^{-1} z_j \sigma_1) z_l = z_l (\sigma_1^{-1} z_j \sigma_1),$
~~$1 \leq  j < l \leq p-1,$

\noindent (R7.2)~$z_l^{-1} (\tau_1^{-1} \xi_j \tau_1) z_l = \tau_1^{-1} \xi_j \tau_1,$
~~$1 \leq  j < l \leq p-1,$

\noindent (R7.3)~$z_j \xi_l = \xi_l z_j,$
~~$1 \leq  j < l \leq p-1,$

\noindent (R8.1)~$(\sigma_1^{-1} z_j \sigma_1^{-1}) z_j = z_j (\sigma_1^{-1} z_j \sigma_1^{-1}),$
~~$1 \leq  j \leq p-1,$

\noindent (R8.2)~$z_j^{-1}( \tau_1^{-1} \xi_j) z_j = \tau_1^{-1} \xi_j,$
~~$1 \leq  j \leq p-1,$

\noindent (R8.3)~$z_j^{-1} \xi_j z_j = \tau_1^{-1} \xi_j \tau_1,$
~~$1 \leq  j \leq p-1.$
\end{prop}

\bigskip

The generators $\sigma_1, \ldots, \sigma_{n-2}, x_1, \ldots, x_{2g}, z_1, \ldots, z_{p-1}$
generate a group isomorphic to $B_{n-1}(\Sigma_{g,p})$ (see also Remark 3.1 from~\cite{Bel}) and 
it is easy to see that the relations (B1), (B2), (R1.1), (R2.1), \ldots, (R8.1)
are a complet set of relations for
$B_{n-1}(\Sigma_{g,p})$. From the other relations we can find the following conjugacy formulae:
\\

\noindent (S1) $\tau_l^{\sigma_k} = \tau_l,$~$k \neq l-1, l,$

\noindent (S2) $\tau_l^{\sigma_{l-1}} = \tau_{l-1},$

\noindent (S3) $\tau_l^{\sigma_l} = \tau_{l+1}^{\tau_l^{-1}},$~$l \neq n-1,$

\noindent (S4) $\tau_i^{x_r} = \tau_i,$~$2 \leq i \leq n-1,$

\noindent (S5) $w_r^{\sigma_i} = w_r,$~$1 \leq r \leq 2g,$~~$1 \leq  i \leq n-2,$

\noindent (S6) $w_r^{x_r} = w_r^{\tau_1},$

\noindent (S7) $\tau_1^{x_r} = \tau_1^{w_r^{-1} \tau_1},$

\noindent (S8) $w_s^{x_r} = w_s^{[w_r^{-1}, \tau_1]},$~$1 \leq s < r \leq 2g,$~$(s, r) \neq (2m-1, 2m),$

\noindent (S9) $w_r^{x_s} = w_r,$~$1 \leq s < r \leq 2g,$~$(s, r) \neq (2m-1, 2m),$

\noindent (S10) $w_{2m-1}^{x_{2m}} = [\tau_1, w_{2m}^{-1}] w_{2m-1} \tau_1,$

\noindent (S11) $w_{2m}^{x_{2m-1}} = \tau_1^{-1} w_{2m},$

\noindent (S12) $\tau_i^{z_j} = \tau_i,$~$2 \leq i \leq n-1,$~$1 \leq j \leq p-1,$

\noindent (S13) $\xi_j^{\sigma_i} = \xi_j,$~$1 \leq i \leq n-2,$~$1 \leq j \leq p-1,$

\noindent (S14) $\xi_j^{x_r} = \xi_j^{[w_r^{-1}, \tau_1]},$~$1 \leq r \leq 2g,$~$1 \leq j \leq p-1,$

\noindent (S15) $w_r^{z_j} = w_r,$~$1 \leq r \leq 2g,$~$1 \leq j \leq p-1,$

\noindent (S16) $\xi_j^{z_l} = \xi_j^{[\xi_l^{-1}, \tau_1]},$~$1 \leq j <l \leq p-1,$

\noindent (S17) $\xi_l^{z_j} = \xi_l,$~$1 \leq j <l \leq p-1,$

\noindent (S18) $\tau_1^{z_j} = [\tau_1, \xi_j^{-1}] \tau_1,$~$1 \leq j \leq p-1,$

\noindent (S19) $\xi_j^{z_j} = \xi_j^{\tau_1},$~$1 \leq j \leq p-1,$

where $a^b = a^{-1} b a$ and $[a, b] = a^{-1} b^{-1} a b.$

\medskip

Let $U_{n-1,g,p}$ be the subgroup of $D_n(\Sigma_{g,p})$ generated by  
$\{\tau_1, \ldots$ $\tau_{n-1}, w_1, \ldots, w_{2g},$ $\xi_1,$ $\ldots,
\xi_{p-1} \}$.

\begin{prop}\label{un} The group $U_{n-1,g,p}$ is a normal subgroup 
of $D_n(\Sigma_{g,p})$ and it is a free group of rank $n+p+2g-2$.
\end{prop}
\begin{proof}
The  statement  was  proven in~\cite[Section.2]{Bel} using the interpretation of 
$U_{n-1,g,p}$ as the fundamental group of the surface $\Sigma_{g,p}$ with $n-1$ points removed. 
\end{proof}

From relations (S1)...(S19) we deduce that
$B_{n-1}(\Sigma_{g,p})$ acts on $U_{n-1,g,p}$ by
conjugacy.

\begin{theorem} \label{theorem3} 
The group $B_{n-1}(\Sigma_{g,p})$ with $n\ge 3$, $g \geq 0$ and $p>0$ acts by conjugacy on the free group
$U_{n-1,g,p}$. Therefore we have a representation $\rho_{U}: B_{n-1}(\Sigma_{g,p}) \to Aut(U_{n-1,g,p})$ 
defined algebraically  as follows:

\noindent -- Generators $\sigma_i,$ $i = 1, \ldots, n-2$:
$$
\sigma_{i} : \left\{
\begin{array}{ll}
\tau_i \longmapsto \tau_{i+1}^{\tau_i^{-1}} \, ; & \\ 
\tau_{i+1} \longmapsto \tau_{i} \, ; & \\ 
\tau_l \longmapsto \tau_{l},   \, l \not= i, i+1, \, ; & \\ 
w_r \longmapsto w_r, \, 1 \leq r \leq 2g \, ; & \\ 
\xi_j \longmapsto \xi_j, \, 1 \leq j \leq p-1.
\end{array} \right.
$$

\noindent -- Generators $x_r,$ $r = 1, \ldots, 2g$:
$$
x_r : \left\{
\begin{array}{ll}
\tau_1  \longmapsto \tau_{1}^{w_r^{-1} \tau_1} \, ; & \\ 
\tau_i   \longmapsto \tau_{i}, \, 2 \leq i \leq n \, ; & \\ 
w_s   \longmapsto  w_{s}^{[w_r^{-1}, \tau_1]},\, s < r, \, (s, r) \neq (2m-1, 2m) \, ; & \\ 
w_{r-1} \longmapsto  [\tau_1, w_r^{-1}] w_{r-1} \tau_1, \; \mbox{if} \; r = 2m \, ; & \\ 
w_r \longmapsto   w_{r}^{\tau_1} \, ; & \\ 
w_s \longmapsto  w_{s},\, r < s,\, (r, s) \neq (2m-1, 2m) \, ; & \\ 
w_{r+1} \longmapsto  \tau_1^{-1} w_{r+1}, \; \mbox{if} \; r = 2m-1 \, ; & \\ 
\xi_j  \longmapsto  \xi_{j}^{[w_r^{-1}, \tau_1]}, \, 1 \leq j \leq p-1\, .
\end{array} \right.
$$

\noindent -- Generators $z_j,$ $j = 1, \ldots, p-1$:
$$
z_j : \left\{
\begin{array}{ll}
\tau_1  \longmapsto  \tau_1^{\xi_j^{-1} \tau_1} \, ; & \\ 
\tau_i \longmapsto  \tau_{i},\, 2 \leq i \leq n \, ; & \\ 
w_r^{z_j}  \longmapsto  w_{r},$~$1 \leq r \leq 2g) \, ; & \\ 
\xi_l^{z_j}  \longmapsto \xi_{l}^{[\xi_j^{-1}, \tau_1]}, \, 1 \leq l < j \leq p-1 \, ; & \\ 
\xi_j^{z_j}    \longmapsto \xi_{j}^{\tau_1} \, ; & \\ 
\xi_l^{z_j}  \longmapsto \xi_{l},\, 1 \leq j < l \leq p-1 \, .
\end{array} \right.
$$
\end{theorem}

In the following we outline a proof of the faithfulness of the representation $\rho_{U}$ of
$B_{n-1}(\Sigma_{g, p})$ given in Theorem ~\ref{theorem3}
 using the interpretation of surface braids as mapping classes.

First, we recall that the mapping class group of a surface $\Sigma_{g,p}$, let us 
denote it by $\mathcal{M}_{g,p}$, is the group of
isotopy classes of orientation-preserving self-homeomorphisms  which
fix the boundary components pointwise.

Let $\mathcal{P}=\{p_1, \dots, p_n\}$ be a set of $n$ distinct points
(\emph{punctures}) in the interior of the surface $\Sigma_{g,p}$.  The punctured
mapping class group of $\Sigma_{g,p}$ relative to $\mathcal{P}$ is
defined  to be the group of isotopy classes of orientation-preserving
self-homeomorphisms which fix the boundary components pointwise, and
which fix $\mathcal{P}$ setwise. This group, that we will denote by $\mathcal{M}_{g,p}^n$, does
not depend on the choice of  $\mathcal{P}$, but just on its cardinal.

We recall also that  a simple closed curve $C$  is \emph{essential} if either it does not bound a
disk or it  bounds a disk containing at least two punctures.

Finally, we denote
$T_C$ as a Dehn twist along a simple closed curve $C$.
Let $C$ and $D$ be two simple closed curves bounding an annulus
containing only the  puncture $p_j$.  We shall say that the
multitwist $T_{C} T_{D}^{-1}$ is a \emph{$j$-bounding pair braid}, also called \emph{spin map}
in~\cite{Bir}.

\vspace*{5pt}

\noindent{\bf Surface braids as mapping classes.} Let $g,p \ge 0$ and let 
$\psi_{n,0}: \mathcal{M}_{g,p}^n \to \mathcal{M}_{g,p}$
be the homomorphism induced by the map which forgets the set
$\mathcal{P}$.
When $p=0$, according to a well-known result of Birman~\cite[Chapter 4.1]{Bir},
the group $B_n(\Sigma_{g})$
is isomorphic to  $\ker \psi_{n,0}$ if $g>1$. 
The statement of Birman's theorem  concerns
  the case of closed surfaces, but the proof  
extends naturally to the case of surfaces with boundary components and the group $B_n(\Sigma_{g,p})$
is isomorphic to  $\ker \psi_{n,0}$ if $g \ge 1$ and $p>0$.

Geometrically the correspondence between  $\ker \psi_{n,0}$ and $B_n(\Sigma_{g,p})$
is realized as follows: given a homeomorphism $h$ of $\Sigma_{g,p}$ isotopic to the identity $Id$ of $\Sigma_{g,p}$ and fixing  boundary components pointwise 
and $\mathcal{P}$ setwise, the track of the punctures $p_1, \ldots, p_n$ under an isotopy from $h$ to $Id$ is the geometric braid corresponding to the homeomorphism $h$.

\begin{theorem}~\label{theorem3annexe}
Let $g \ge 0$, $p>0$ and $n\ge 3$. 
The representation $\rho_{U}: B_{n-1}(\Sigma_{g,p}) \to {\rm Aut} (U_{n-1,g,p})$ is faithful.
\end{theorem}
\begin{proof}
First we prove the claim when  $g$ is greater or equal then $1$.
According to Birman's result we can represent the group $D_{n}(\Sigma_{g,p})$ as a  (normal)
subgroup  of $\mathcal{M}_{g,p}^n$, more precisely as the subgroup of mapping classes in  $\ker \psi_{n,0}$
sending the $n$th puncture into itself (see also Remark~18 from~\cite{BGG}). In particular generators 
$\tau_1, \ldots, \tau_{n-1},~w_1, \ldots, w_{2g},~\xi_1, \ldots,
\xi_{p-1},
$
of $D_{n}(\Sigma_{g,p})$ (see Proposition~\ref{propdn1}) correspond to $n$-bounding pair braids in $\mathcal{M}_{g,p}^n$.

Let $\psi_{n,n-1}: \mathcal{M}_{g,p}^n \to \mathcal{M}_{g,p}^{n-1}$
be the homomorphism  forgetting the last puncture.
Now, let us consider the exact sequence associate to the restriction of $\psi_{n,n-1}$
to $D_{n}(\Sigma_{g,p})$. The image of $D_{n}(\Sigma_{g,p})$ by  $\psi_{n,n-1}$ coincides with $\ker \psi_{n-1,0}$
which is isomorphic to $B_{n-1}(\Sigma_{g,p})$. On the other hand, $\ker \psi_{n,n-1} \cap D_{n}(\Sigma_{g,p})$
 is the subgroup of 
 $\mathcal{M}_{g,p}^n$ generated by $n$-bounding pair braids 
$\tau_1, \ldots, \tau_{n-1},~w_1, \ldots, w_{2g},~\xi_1, \ldots,
\xi_{p-1}$ and therefore it is isomorphic to  $U_{n-1,g,p}$. 

The group $B_{n-1}(\Sigma_{g,p})$ embeds naturally
in $B_{n}(\Sigma_{g,p})$ by sending generators of $B_{n-1}(\Sigma_{g,p})$ in corresponding ones
of $B_{n}(\Sigma_{g,p})$ and therefore $B_{n-1}(\Sigma_{g,p})$ can be considered also as a subgroup
of $\mathcal{M}_{g,p}^n$. Moreover, the action considered in Theorem~\ref{theorem3} corresponds to the action by
conjugacy of $B_{n-1}(\Sigma_{g,p})$, seen as a subgroup of $\mathcal{M}_{g,p}^n$, on $U_{n-1,g,p}$.
Since $U_{n-1,g,p}$ is a normal subgroup of $\mathcal{M}_{g,p}^n$
we can define a map $\Theta: {\rm Inn}(\mathcal{M}_{g,p}^n) \to {\rm Aut}(U_{n-1,g,p})$.
We prove that $\Theta$ is injective and therefore that the action by
conjugacy of $B_{n-1}(\Sigma_{g,p})$ on $U_{n-1,g,p}$ is faithful.
Let $g \in  {\rm Inn}(\mathcal{M}_{g,p}^n)$ such that $\Theta(g)=1$ and $C$ be an essential curve. 
We can associate to the curve $C$
another simple closed curve
$D$ such that they bound an annulus containing only the  puncture $p_n$.
The mapping class  $T_C T_{D}^{-1}$ is then a $n$-bounding pair
braid.

Since $\Theta(g)=1$, then $g \, T_C T_{D}^{-1} g^{-1}= T_C T_{D}^{-1}$ and from a simple argument on the
index of intersection of curves   (see for instance Proposition~2.10 of~\cite{belsph})
one can easily deduce that
$g(C)=C$. Since $g(C)=C$ for any essential curve $C$
it follows that $g$ is isotopic to the identity
(Lemma 5.1 and Theorem 5.3 of \cite{IM}). Therefore $\Theta$ is injective and in particular the representation defined in Theorem~\ref{theorem3} is faithful when $g\ge 1$.

The only case left is when the genus is equal to zero.
We recall that, for $p>0$, the group $B_{n}(\Sigma_{0,p})$ is isomorphic 
to the subgroup $B_{n+p-1, p-1}$ of $B_{n+p-1}$
fixing the last $p-1$ strands. In this case our representation coincides with
Artin representation of $B_{n+p-1}$ in ${\rm Aut}(F_{n+p-1})$
restricted to the subgroup $B_{n+p-1, p-1}$.
\end{proof}

\begin{cor}
The group $B_m(\Sigma_{g,p})$ 
is residually finite for $m\geq 1$, $g \geq 0$ and
$p > 0$.
\end{cor}
\begin{proof}
Baumslag and Smirnov
 proved (see for instance \cite[Theorem 4.8]{LS})  that any finitely presented group
 which is isomorphic to a subgroup of automorphisms of a free group of finite rank is residually finite.
 Therefore for $m>1$ the statement is a Corollary of Theorem~\ref{theorem3annexe}.
 In the case $m=1$, the group  $B_1(\Sigma_{g,p})$ is isomorphic to the fundamental group
 of  $\Sigma_{g,p}$ which is free and therefore residually finite.
\end{proof}

Remark also that from Theorem~\ref{theorem3annexe} one can derive, using Fox derivatives,
a Burau representation for $B_n(\Sigma_{g,p})$. Since the restriction on $B_n$ of such representation coincides with the usual Burau representation of $B_n$, the Burau representation for $B_n(\Sigma_{g,p})$ obtained from  Theorem~\ref{theorem3annexe} is not faithful for $n\ge 5$.  


\section{Surface braids as outer automorphisms of free groups}

As recalled in the Introduction, any element  $\beta$ of $B_n \subset {\rm Aut}(F_n)$ fixes the product
$x_1 x_2 \ldots x_n$ of generators of $F_n$. We prove a similar statement  for the
group $B_n(\Sigma_{g,p})$, for $p > 0$.

\begin{prop} \label{prop1} Let $p>0$ and let $U_{n-1,g,p}$ be the free group of rank  $n+p+2g-2$
defined above. 
Any element $\beta$ in  $B_{n-1}(\Sigma_{g,p}) \subset {\rm Aut}(U_{n-1,g,p})$ fixes the
product
$$
A = \tau_{n-1}^{-1} \ldots \tau_{2}^{-1} \tau_1^{-1} \xi_1 \xi_2 \ldots \xi_{p-1} [w_1^{-1}, w_2] \ldots
[w_{2g - 1}^{-1}, w_{2g}].
$$
\end{prop}
\begin{proof} In order  to prove the claim it suffices to verify that any generator  of
$
B_{n-1}(\Sigma_{g,p})$,
considered as an automorphism of $U_{n-1,g,p}$, fixes  the element $A$.
\\

\noindent {\it Case 1:} generators $\sigma_i$, $1 \leq i \leq n-2$.
The group $B_{n-1}$ is a subgroup of $B_{n-1}(\Sigma_{g,p})$ and by Artin's theorem $B_{n-1}$ is a
subgroup of ${\rm Aut}(F_{n-1}),$ $F_{n-1} = \langle \tau_1, \tau_2, \ldots, \tau_{n-1} \rangle,$ and fixes
the product $\tau_1 \tau_2 \ldots \tau_{n-1}$. Hence, the generator $\sigma_i$  also fixes the product
$\tau_{n-1}^{-1} \ldots \tau_{2}^{-1} \tau_1^{-1}$.

From Theorem~\ref{theorem3} it follows that
$$
w_r^{\sigma_i} = w_r, 1 \leq r \leq 2g;~~\xi_j^{\sigma_i} = \xi_j, 1 \leq j \leq p-1.
$$
Hence, $A^{\sigma_i} = A.$
\\

\noindent {\it Case 2:} generators $x_r$, $1 \leq r \leq 2g$.
By Theorem~\ref{theorem3}  we have
$$
(\tau_{n-1}^{-1} \ldots \tau_{2}^{-1} \tau_1^{-1})^{x_r} = \tau_{n-1}^{-1} \ldots \tau_{2}^{-1}
(\tau_1^{-1} w_r \tau_1^{-1} w_r^{-1} \tau_1)
$$
and
$$
(\xi_1 \ldots \xi_{p-1})^{x_r} = (\xi_1 \ldots \xi_{p-1})^{[w_r^{-1}, \tau_1]}.
$$

Let us consider
$$
([w_1^{-1}, w_2] \ldots [w_{2g - 1}^{-1}, w_{2g}])^{x_r}.
$$
We will distinguish two cases: $r=2m-1$ and $r=2m,$ where $m=1, \ldots g.$
In the first case we have the following equalities:
$$
([w_1^{-1}, w_2] \ldots [w_{2m - 3}^{-1}, w_{2m-2}])^{x_{2m-1}} =
([w_1^{-1}, w_2] \ldots [w_{2m - 3}^{-1}, w_{2m-2}])^{[w_{2m-1}^{-1}, \tau_1]},
$$

$$
[w_{2m - 1}^{-1}, w_{2m}]^{x_{2m-1}} = [w_{2m-1}^{-\tau_1}, \tau_1^{-1} w_{2m}] =
[\tau_1, w_{2m-1}^{-1}] [w_{2m - 1}^{-1}, w_{2m}],
$$

$$
([w_{2m+1}^{-1}, w_{2m+2}] \ldots [w_{2g - 1}^{-1}, w_{2g}])^{x_{2m-1}} =
[w_{2m+1}^{-1}, w_{2m+2}] \ldots [w_{2g - 1}^{-1}, w_{2g}].
$$

In the second case we have the following equalities:

$$
([w_1^{-1}, w_2] \ldots [w_{2m - 3}^{-1}, w_{2m-2}])^{x_{2m}} =
([w_1^{-1}, w_2] \ldots [w_{2m - 3}^{-1}, w_{2m-2}])^{[w_{2m}^{-1}, \tau_1]},
$$

$$
[w_{2m - 1}^{-1}, w_{2m}]^{x_{2m}} = [\tau^{-1}_1 w_{2m-1}^{-1} [w_{2m}^{-1}, \tau_1],
w_{2m}^{\tau_1}] =
[\tau_1, w_{2m}^{-1}] [w_{2m - 1}^{-1}, w_{2m}],
$$

$$
([w_{2m+1}^{-1}, w_{2m+2}] \ldots [w_{2g - 1}^{-1}, w_{2g}])^{x_{2m}} =
[w_{2m+1}^{-1}, w_{2m+2}] \ldots [w_{2g - 1}^{-1}, w_{2g}].
$$

In  both  cases we have
$$
([w_1^{-1}, w_2] \ldots [w_{2g - 1}^{-1}, w_{2g}])^{x_r} =
[\tau_1, w_r^{-1}] [w_1^{-1}, w_2] \ldots [w_{2g - 1}^{-1}, w_{2g}].
$$
Hence, the element $x_r$ acts by conjugation on $A$ as follows
$$
A^{x_r} = \tau_{n-1}^{-1} \ldots \tau_{2}^{-1}
\tau_1^{-1} [w_r^{-1}, \tau_1]
(\xi_1 \ldots \xi_{p-1})^{[w_r^{-1}, \tau_1]}
 [\tau_1, w_r^{-1}] [w_1^{-1}, w_2] \ldots $$
$[w_{2g - 1}^{-1}, w_{2g}] ,
$

and then it is easy to check that $A^{x_r} = A.$
\\

\noindent {\it Case 3:} generators $z_j$, $1 \leq j \leq p-1$.
By Theorem~\ref{theorem3}  we have

$$
(\tau_{n-1}^{-1} \ldots \tau_{2}^{-1} \tau_1^{-1})^{z_j} = \tau_{n-1}^{-1} \ldots \tau_{2}^{-1}
\tau_1^{-1} [\xi_j^{-1}, \tau_1],
$$

$$
(\xi_1 \ldots \xi_{p-1})^{z_j} = (\xi_1 \ldots \xi_{j-1})^{[\xi_j^{-1}, \tau_1]}
\xi_j^{\tau_1} \xi_{j+1} \ldots \xi_{p-1} = [\xi_j^{-1}, \tau_1]^{-1} \xi_{1} \ldots \xi_{p-1},
$$

and

$$
([w_1^{-1}, w_2] \ldots [w_{2g - 1}^{-1}, w_{2g}])^{z_j} =
 [w_1^{-1}, w_2] \ldots [w_{2g - 1}^{-1}, w_{2g}].
$$
Hence, $A^{z_j} = A$  and the Proposition follows.
\end{proof}

\vspace{5pt}
We recall that $\rho_U: B_{n-1}(\Sigma_{g,p}) \to {\rm Aut} (U_{n-1,g,p})$ is the representation
of $B_{n-1}(\Sigma_{g,p})$ defined in Theorem~\ref{theorem3} and let 
$p: {\rm Aut} (U_{n-1,g,p})\to {\rm Out} (U_{n-1,g,p})$ be the canonical projection.  

\begin{theorem}~\label{out}
 The representation  $p\circ \rho_U: B_{n-1}(\Sigma_{g,p}) \to {\rm Out} (U_{n-1,g,p})$  is faithful
for $p>0$ when $g>0$ and for $p>2$ when $g=0$.
\end{theorem}
\begin{proof}
Since the representation $\rho_U: B_{n-1}(\Sigma_{g,p}) \to {\rm Aut} (U_{n-1,g,p})$ 
is faithful (Theorem~\ref{theorem3annexe}),  we can identify $B_{n-1}(\Sigma_{g,p})$
with  $\rho_U(B_{n-1}(\Sigma_{g,p}))$.

Now suppose that there exists  $\beta \in B_{n-1}(\Sigma_{g,p}) \cap {\rm Inn} (U_{n-1,g,p})$.
From Proposition~\ref{prop1} one deduces that $\beta$ is a conjugation by a power $m$ of 
$$
A = \tau_{n-1}^{-1} \ldots \tau_{2}^{-1} \tau_1^{-1} \xi_1 \xi_2 \ldots \xi_{p-1} [w_1^{-1}, w_2] \ldots
[w_{2g - 1}^{-1}, w_{2g}].
$$
Now, let $g$ be a generator of $B_{n-1}(\Sigma_{g,p})$.
Since all elements of $B_{n-1}(\Sigma_{g,p})$ fix the element $A$  we  deduce the following
equalities:
$$
g^{-1}( \beta (g (x))) = g^{-1}( A^m g (x) A^{-m})= A^m x A^{-m}= \beta(x),
$$
for any $x \in U_{n-1,g,p}$.
One deduces that $g^{-1} \beta g= \beta$  for any generator $g$
of $B_{n-1}(\Sigma_{g,p})$ and  therefore $\beta$ belongs to the center of $B_{n-1}(\Sigma_{g,p})$.
Since   $B_{n-1}(\Sigma_{g,p})$
(with $n\ge 2$) has trivial center for $p>0$ when $g>0$~\cite{RP} and for $p>2$ when $g=0$, the 
intersection $B_{n-1}(\Sigma_{g,p}) \cap {\rm Inn} (U_{n-1,g,p})$ is trivial and the claim follows.
\end{proof}

\begin{rem} \label{outrem}
 Let $\phi: B_n \to {\rm Out} (F_n)$ the representation obtained composing  Artin representation
of $B_n$ in ${\rm Aut} (F_n)$  with the canonical projection of ${\rm Aut} (F_n)$  in ${\rm Out} (F_n)$.
Such representation is not faithful and it is easy to see that the kernel is the center of $B_n$.
\end{rem}


\section{Surface braid groups of non-orientable surfaces with boundary components}
Let $N_{g,p}$ be a non-orientable surface of genus $g \geq 1,$
with $p>0$ boundary components.
Let $\sigma_1, \dots, \sigma_{n-1}$ be the usual generators of $B_n$ and
$a_1, \ldots a_{g}, z_1, \ldots z_{p-1}$ be the usual generators of the fundamental group of
$N_{g,p}$. As in previous section we can consider $\sigma_1, \dots, \sigma_{n-1}$  and
$a_1, \ldots a_{g}, z_1, \ldots z_{p-1}$ as elements of  $B_n(N_{g, p})$.
According to~\cite{Bel} the group $B_n(N_{g, p})$ admits a presentation with generators:
$$
\sigma_1, \ldots \sigma_{n-1}, a_1, \ldots a_{g}, z_1, \ldots z_{p-1},
$$
and relations:

\noindent -- Braid relations:

$\sigma_i \sigma_{i+1} \sigma_i = \sigma_{i+1} \sigma_i \sigma_{i+1},~~1 \leq i \leq n-2,$

$\sigma_i \sigma_{j} = \sigma_{j} \sigma_i,~~|i - j| > 1,~~1 \leq i,j \leq n-1,$

\noindent -- Mixed relations:

\noindent (R1)~$a_r \sigma_i = \sigma_i a_r,$~~$i \neq 1,$~~$1 \leq r \leq g,$

\noindent (R2)~$\sigma_1^{-1} a_r \sigma_1^{-1} a_r = a_r \sigma_1^{-1} a_r \sigma_1,$
~~$1 \leq r \leq g,$

\noindent (R3)~$(\sigma_1^{-1} a_s \sigma_1) a_r = a_r (\sigma_1^{-1} a_s \sigma_1)$,
$1 \leq s < r \leq g,$

\noindent (R4)~$z_j \sigma_i = \sigma_i z_j,$~~$i \neq 1,$~~$1 \leq j \leq p-1,$

\noindent (R5)~$(\sigma_1^{-1} z_i \sigma_1) a_r = a_r (\sigma_1^{-1} z_i \sigma_1)$, $1 \leq  r \leq g$, $1 \leq i \leq p-1$, $n > 1,$

\noindent (R6)~$(\sigma_1^{-1} z_j \sigma_1) z_l = z_l (\sigma_1^{-1} z_j \sigma_1),$
~~$1 \leq  j < l \leq p-1,$

\noindent (R7)~$(\sigma_1^{-1} z_j \sigma_1^{-1}) z_j = z_j (\sigma_1^{-1}
z_j \sigma_1^{-1})$, $1 \leq  j \leq p-1.$
\\

As in previous section let us consider the natural projection of 
$\pi : B_n(N_{g,p}) \longrightarrow S_n$,
which associates to any braid the corresponding permutation.
This projection map $\sigma_i$ in the corresponding transposition and 
generators $a_1, \ldots, a_g, z_1, \ldots, z_{p-1}$ into the identity. 
As before, let  $D_n(N_{g,p}) = \pi^{-1}(S_{n-1})$
and let $m_l = \sigma_{n-1} \ldots$ $ \sigma_l,~l = 1, \ldots, n-1$, $m_n=1$.
The set $M_n = \{ m_l ~|~1 \leq l \leq n \}$ is a Schreier set of coset representatives of
$D_n(N_{g,p})$ in $B_n(N_{g,p})$ and the group $D_n(N_{g,p})$ is generated by
$$
s_{\lambda, a} = \lambda a \cdot (\overline{\lambda a})^{-1},
$$
where $\lambda$ runs over the set $M_n$ and $a$ runs over the set of generators of
$B_n(N_{g,p})$.
\\

\noindent {\bf Case 1.} If $a \in \{ \sigma_1, \ldots, \sigma_{n-1} \}$, then we find the generators
$$
\tau_k = \sigma_{n-1} \ldots \sigma_{k+1} \sigma_k^2 \sigma_{k+1}^{-1} \ldots \sigma_{n-1}^{-1},~
k=1,\ldots,n-2,~\tau_{n-1}=\sigma_{n-1}^2.
$$
\\

\noindent {\bf Case 2.} If $a \in \{ a_1, \ldots, a_{g} \}$, then we find the generators
$$
w_r = \sigma_{n-1} \ldots \sigma_{1} a_r \sigma_{1}^{-1} \ldots \sigma_{n-1}^{-1},~
r=1,\ldots,g.
$$
\\

\noindent {\bf Case 3.} If $a \in \{ z_1, \ldots, z_{p-1} \}$, then we find the generators
$$
\xi_j = \sigma_{n-1} \ldots \sigma_{1} z_j \sigma_{1}^{-1} \ldots \sigma_{n-1}^{-1},~
j=1,\ldots,p-1.
$$

Using the same argument as in the
orientable case one can find the following group presentation for $D_n(N_{g, p})$.

\begin{prop} \label{theorem4}
The group $D_n(N_{g, p}),$ $n \geq 1$ admits a presentation with the generators
$$
\sigma_1, \ldots, \sigma_{n-2}, a_1, \ldots, a_{g}, z_1, \ldots, z_{p-1},
\tau_1, \ldots, \tau_{n-1},~w_1, \ldots, w_{g},~\xi_1, \ldots, \xi_{p-1},
$$
and   the following relations:

\noindent -- Braid relations:
\\

\noindent (B1)~$\sigma_i \sigma_{i+1} \sigma_i = \sigma_{i+1} \sigma_i \sigma_{i+1},~~1 \leq i \leq n-3,$

\noindent (B2)~$\sigma_i \sigma_{j} = \sigma_{j} \sigma_i,~~|i - j| > 1,~~1 \leq i,j \leq n-2,$

\noindent (B3)~$\sigma_k^{-1} \tau_l \sigma_k = \tau_l,~k \neq l-1, l,$

\noindent (B4)~$\sigma_{l-1}^{-1} \tau_l \sigma_{l-1} = \tau_{l-1},$

\noindent (B5)~$\sigma_l^{-1} \tau_l \sigma_l = \tau_l \tau_{l+1} \tau_l^{-1},~l \neq n-1.$

\noindent -- Mixed relations:
\\

\noindent (R1.1)~$a_r \sigma_i = \sigma_i a_r,$~~$2 \leq i \leq n-2,$~~$1 \leq r \leq g,$

\noindent (R1.2)~$a_r \tau_i = \tau_i a_r,$~~$2 \leq i \leq n-1,$

\noindent (R1.3)~$w_r \sigma_i = \sigma_i w_r,$~~$1 \leq r \leq g,$~~$1 \leq i \leq n-2,$

\noindent (R2.1)~$(\sigma_1^{-1} a_r \sigma_1^{-1}) a_r = a_r (\sigma_1^{-1} a_r \sigma_1),$
~~$1 \leq r \leq g,$

\noindent (R2.2)~$a_r^{-1} (\tau_1^{-1} w_r) a_r = \tau_1^{-1} w_r \tau_1,$
~~$1 \leq r \leq g,$

\noindent (R2.3)~$a_r^{-1} w_r a_r = \tau_1^{-1} w_r,$
~~$1 \leq r \leq g,$

\noindent (R3.1)~$(\sigma_1^{-1} a_s \sigma_1) a_r = a_r (\sigma_1^{-1} a_s \sigma_1),$
~~$1 \leq s < r \leq g,$

\noindent (R3.2)~$a_r^{-1} (\tau_1^{-1} w_s \tau_1) a_r = \tau_1^{-1} w_s \tau_1,$
~~$1 \leq s < r \leq g,$

\noindent (R3.3)~$a_s w_r = w_r a_s,$
~~$1 \leq s < r \leq g,$

\noindent (R4.1)~$z_j \sigma_i = \sigma_i z_j,$~~$2 \leq i \leq n-2,$~~$1 \leq j \leq p-1,$

\noindent (R4.2)~$z_j \tau_i = \tau_i z_j,$~~$2 \leq i \leq n-1,$~~$1 \leq j \leq p-1,$

\noindent (R4.3)~$\xi_j \sigma_i = \sigma_i \xi_j,$~~$1 \leq j \leq p-1,$~~$1 \leq i \leq n-2,$

\noindent (R5.1)~$(\sigma_1^{-1} z_j \sigma_1) a_r = a_r (\sigma_1^{-1} z_j \sigma_1),$
~~$1 \leq  r \leq g,$~~$1 \leq j \leq p-1,$

\noindent (R5.2)~$a_r^{-1} (\tau_1^{-1} \xi_j \tau_1) a_r = \tau_1^{-1} \xi_j \tau_1,$
~~$1 \leq  r \leq g,$~~$1 \leq j \leq p-1,$

\noindent (R5.3)~$z_j w_r = w_r z_j,$
~~$1 \leq  r \leq g,$~~$1 \leq j \leq p-1,$

\noindent (R6.1)~$(\sigma_1^{-1} z_j \sigma_1) z_l = z_l (\sigma_1^{-1} z_j \sigma_1),$
~~$1 \leq  j < l \leq p-1,$

\noindent (R6.2)~$z_l^{-1} (\tau_1^{-1} \xi_j \tau_1) z_l = \tau_1^{-1} \xi_j \tau_1,$
~~$1 \leq  j < l \leq p-1,$

\noindent (R6.3)~$z_j \xi_l = \xi_l z_j,$
~~$1 \leq  j < l \leq p-1,$

\noindent (R7.1)~$(\sigma_1^{-1} z_j \sigma_1^{-1}) z_j = z_j (\sigma_1^{-1} z_j \sigma_1^{-1}),$
~~$1 \leq  j \leq p-1,$

\noindent (R7.2)~$z_j^{-1}( \tau_1^{-1} \xi_j) z_j = \tau_1^{-1} \xi_j,$
~~$1 \leq  j \leq p-1,$

\noindent (R7.3)~$z_j^{-1} \xi_j z_j = \tau_1^{-1} \xi_j \tau_1,$
~~$1 \leq  j \leq p-1.$
\end{prop}

It is easy to see that the relations (B1), (B2), (R1.1), (R2.1), \ldots, (R7.1)
are defining relations of
$B_{n-1}(N_{g,p})$. From the other relations we can find the following conjugacy formulae:

\noindent (S1) $\tau_l^{\sigma_k} = \tau_l,$~$k \neq l-1, l,$

\noindent (S2) $\tau_l^{\sigma_{l-1}} = \tau_{l-1},$

\noindent (S3) $\tau_l^{\sigma_l} = \tau_{l+1}^{\tau_l^{-1}},$~$l \neq n-1,$

\noindent (S4) $\tau_i^{a_r} = \tau_i,$~$2 \leq i \leq n-1,$

\noindent (S5) $w_r^{\sigma_i} = w_r,$~$1 \leq r \leq g,$~~$1 \leq  i \leq n-2,$

\noindent (S6) $\tau_1^{a_r} = (\tau_1^{-1})^{w_r^{-1} \tau_1},$

\noindent (S7) $w_r^{a_r} = \tau_1^{-1} w_r,$~$1 \leq r \leq g,$

\noindent (S8) $w_s^{a_r} = w_s^{w_r \tau_1 w_r^{-1} \tau_1},$~$1 \leq s < r \leq g,$

\noindent (S9) $w_r^{a_s} = w_r,$~$1 \leq s < r \leq g,$

\noindent (S10) $\tau_i^{z_j} = \tau_i,$~$2 \leq i \leq n-1,$~$1 \leq j \leq p-1,$

\noindent (S11) $\xi_j^{\sigma_i} = \xi_j,$~$1 \leq i \leq n-2,$~$1 \leq j \leq p-1,$

\noindent (S12) $\xi_j^{a_r} = \xi_j^{w_r \tau_1 w_r^{-1} \tau_1},$~$1 \leq r \leq g,$
~$1 \leq j \leq p-1,$

\noindent (S13) $w_r^{z_j} = w_r,$~$1 \leq r \leq g,$~$1 \leq j \leq p-1,$

\noindent (S14) $\xi_j^{z_l} = \xi_j^{[\xi_l^{-1}, \tau_1]},$~$1 \leq j <l \leq p-1,$

\noindent (S15) $\xi_l^{z_j} = \xi_l,$~$1 \leq j <l \leq p-1,$

\noindent (S16) $\tau_1^{z_j} = [\tau_1, \xi_j^{-1}] \tau_1,$~$1 \leq j \leq p-1,$

\noindent (S17) $\xi_j^{z_j} = \xi_j^{\tau_1},$~$1 \leq j \leq p-1.$
\\

Consider the subgroup
$$
W_{n,g,p} = \langle \tau_1, \ldots, \tau_{n-1}, w_1, \ldots, w_{g}, \xi_1, \ldots,
\xi_{p-1} \rangle
$$
of $D_{n+1}(N_{g,p})$. The group $W_{n,g,p}$   is free (see \cite{Bel}) and from  (S1)\noindent --(S17) we see that $B_{n}(N_{p,g})$ acts on
$W_{n,g,p}$ by
conjugacy.

\begin{theorem} \label{theorem5a}
The group $B_n(N_{g, p})$, for  $n>1$, $p \geq 1$ and $g\ge 1$ acts by conjugation on the free group
$W_{n,g,p}$
and the induced representation $\rho_W: B_n(N_{g, p}) \to {\rm Aut}(W_{n,g,p})$ is defined as follows:

\noindent -- Generators $\sigma_i,$ $i = 1, \ldots, n-1$:
$$
\sigma_i : \left\{
\begin{array}{ll}
\tau_i  \longmapsto \tau_{i+1}^{\tau_i^{-1}}\, ; & \\ 
\tau_{i+1}  \longmapsto \tau_{i}\, ; & \\ 
\tau_l \longmapsto  \tau_{l},$~$l \neq i, i+1\, ; & \\ 
w_r  \longmapsto w_r,\, 1 \leq r \leq g\, ; & \\ 
\xi_j  \longmapsto \xi_j,$~$1 \leq j \leq p-1\, .
\end{array} \right.
$$

\noindent -- Generators $a_r,$ $r = 1, \ldots, g$:
$$
a_r : \left\{
\begin{array}{ll}
\tau_1   \longmapsto (\tau_{1}^{-1})^{w_r^{-1} \tau_1}\, ; & \\ 
\tau_i   \longmapsto \tau_{i},\, 2 \leq i \leq n\, ; & \\ 
w_s  \longmapsto  w_{s}^{w_r \tau_1 w_r^{-1} \tau_1},$~$1 \leq s < r \leq g\, ; & \\ 
w_r  \longmapsto  \tau_1^{-1} w_{r}\, ; & \\ 
w_s   \longmapsto w_{s},$~$ 1 \leq r < s \leq g\, ; & \\ 
\xi_j   \longmapsto \xi_{j}^{w_r \tau_1 w_r^{-1} \tau_1},$~$1 \leq j \leq p-1.
\end{array} \right.
$$

\noindent -- Generators  $z_j,$ $j = 1, \ldots, p-1$:
$$
z_j : \left\{
\begin{array}{ll}
\tau_1   \longmapsto [\tau_{1}, \xi_j^{-1}] \tau_1\, ; & \\ 
\tau_i   \longmapsto \tau_{i},$~$2 \leq i \leq n\, ; & \\ 
w_r   \longmapsto w_{r},\, 1 \leq r \leq g\, ; & \\ 
\xi_l  \longmapsto \xi_{l}^{[\xi_j^{-1}, \tau_1]},\, 1 \leq l < j \leq p-1\, ; & \\ 
\xi_j   \longmapsto \xi_{j}^{\tau_1}\, ; & \\ 
\xi_l  \longmapsto  \xi_{l},\, 1 \leq j < l \leq p-1.
\end{array} \right.
$$
\end{theorem}

We don't know if the  representation  $\rho_W: B_n(N_{g, p}) \to {\rm Aut}(W_{n,g,p})$
is faithful or not. 

\begin{prop} \label{prop2}
Let $\beta$ be in $B_n(N_{g, p})$.
The  element  $\rho_W(\beta)$ in ${\rm Aut}(W_{n,g,p})$ fixes the
product
$$
A = \tau_{n-1}^{-1} \ldots \tau_{2}^{-1} \, \tau_1^{-1} \, \xi_1 \, \xi_2 \ldots \xi_{p-1} \,
w_1^{2}  \, w_2^2 \ldots
 w_{g}^2.
$$
\end{prop}
\begin{proof} As above it is enough to prove that images of generators of $B_n(N_{g, p})$
fix $A$.

\noindent {\it Case 1:} generators $\sigma_i,$ $1 \leq i \leq n-1$.
The group $B_n$ is a subgroup of $B_n(N_{g,p})$ and by Artin's theorem $B_n$ is a
subgroup of ${\rm Aut}(F_n)$, $F_n = \langle \tau_1, \tau_2, \ldots, \tau_{n-1} \rangle$ and fixes
the product $\tau_1 \tau_2 \ldots \tau_{n-1}$. Hence, $\rho_W(\sigma_i)$  also fixes the product
$\tau_{n-1}^{-1} \ldots \tau_{2}^{-1} \tau_1^{-1}$.

From Theorem~\ref{theorem5a} it follows that
$$
w_r^{\sigma_i} = w_r, 1 \leq r \leq g;~~\xi_j^{\sigma_i} = \xi_j, 1 \leq j \leq p-1.
$$
Hence, $A^{\sigma_i} = A.$
\\

\noindent {\it Case 2:} generators $a_r,$ $1 \leq r \leq g$. By Theorem~\ref{theorem5a} we have that
$$
(\tau_{n-1}^{-1} \ldots \tau_{2}^{-1} \tau_1^{-1})^{a_r} = \tau_{n-1}^{-1} \ldots \tau_{2}^{-1}
(\tau_1^{w_r^{-1} \tau_1}),
$$
and
$$
(\xi_1 \ldots \xi_{p-1})^{a_r} = (\xi_1 \ldots \xi_{p-1})^{w_r \tau_1 w_r^{-1} \tau_1},
$$
and
$$
(w_1^2 \ldots w_{g}^2)^{a_r} = (w_1^2 \ldots w_{r-1}^2)^{w_r \tau_1 w_r^{-1} \tau_1}
(\tau_1^{-1} w_r \tau_1^{-1} w_r)  (w_{r+1}^2 \ldots w_{g}^2).
$$

Hence, the element $a_r$ acts on $A$ as follows:
$$
A^{a_r} = \tau_{n-1}^{-1} \ldots \tau_{2}^{-1} \,
\tau_1^{-1} \,
\xi_1  \ldots \xi_{p-1} \, w_1^{2} \, w_2^2 \ldots  w_{g}^2,
$$
and  $A^{a_r} = A.$

\noindent {\it Case 3.} Generators $z_j,$ $1 \leq j \leq p-1$. By Theorem~\ref{theorem5a} we deduce that

$$
(\tau_{n-1}^{-1} \ldots \tau_{2}^{-1} \tau_1^{-1})^{z_j} = \tau_{n-1}^{-1} \ldots \tau_{2}^{-1}
\, (\tau_1^{-1} \, \xi_j \, \tau_1^{-1} \, \xi_j^{-1} \, \tau_1),
$$

$$
(\xi_1 \ldots \xi_{p-1})^{z_j} = (\xi_1 \ldots \xi_{j-1})^{\xi_j \tau_1^{-1} \xi_j^{-1} \tau_1}
\, (\xi_j^{\tau_1}) \, \xi_{j+1} \ldots \xi_{p-1},
$$
and
$$
(w_1^{2} \, w_2^2 \ldots w_{g}^2)^{z_j} =
w_1^{2} \, w_2^2 \ldots w_{g}^2.
$$
Hence, $A^{z_j} = A$.
\end{proof}


\section{Braid groups of closed surfaces}

Let
$\Sigma_g = \Sigma_{g, 0}$
be a closed orientable surface of genus $g \geq 1$.
The group $B_n(\Sigma_{g}),$ $n>1$ admits a group presentation with  generators:
$$
\sigma_1, \ldots \sigma_{n-1}, x_1, \ldots x_{2g},
$$
and   relations:

\noindent -- Braid relations:

$\sigma_i \sigma_{i+1} \sigma_i = \sigma_{i+1} \sigma_i \sigma_{i+1},~~1 \leq i \leq n-2,$

$\sigma_i \sigma_{j} = \sigma_{j} \sigma_i,~~|i - j| > 1,~~1 \leq i,j \leq n-1,$

\noindent -- Mixed relations:

\noindent (R1)~$x_r \sigma_i = \sigma_i x_r,$~~$i \neq 1,$~~$1 \leq r \leq 2g,$

\noindent (R2)~$(\sigma_1^{-1} x_r \sigma_1^{-1}) x_r = x_r (\sigma_1^{-1} x_r \sigma_1^{-1}),$
~~$1 \leq r \leq 2g,$

\noindent (R3)~$(\sigma_1^{-1} x_s \sigma_1) x_r = x_r (\sigma_1^{-1} x_s \sigma_1),$
~~$1 \leq s < r \leq 2g,$~~$(s, r) \neq (2m-1, 2m),$

\noindent (R4)~$(\sigma_1^{-1} x_{2m-1} \sigma_1^{-1}) x_{2m} = x_{2m} (\sigma_1^{-1} x_{2m-1}
\sigma_1),$
~~$1 \leq m \leq g,$

\noindent (TR5)~$[x_1^{-1}, x_2] [x_3^{-1}, x_4] \ldots [x_{2g-1}^{-1}, x_{2g}] = \sigma_1
\sigma_2 \ldots \sigma_{n-1}^2 \ldots \sigma_2 \sigma_1.$
\\

As before, let $\pi$  be the natural projection of $B_n(\Sigma_{g})$
in $S_n$,
let $D_n(\Sigma_{g}) = \pi^{-1}(S_{n-1})$
 and
$m_l = \sigma_{n-1} $ $\ldots \sigma_l,~l = 1, $ $\ldots, n-1,~m_n=1$.
The set $M_n = \{ m_l ~|~1 \leq l \leq n \}$ is a Schreier set of coset representatives of
$D_n(\Sigma_{g})$ in $B_n(\Sigma_{g})$ and $D_n(\Sigma_{g})$ is generated by
$$
\sigma_1, \sigma_2, \ldots , \sigma_{n-2}, \tau_1, \tau_2, \ldots , \tau_{n-1},
x_1, x_2, \ldots , x_{2g}, w_1, w_2, \ldots , w_{2g},
$$
where
$$
\tau_k = \sigma_{n-1} \ldots \sigma_{k+1} \sigma_k^2 \sigma_{k+1}^{-1} \ldots \sigma_{n-1}^{-1},~
k=1,\ldots,n-2,~\tau_{n-1}=\sigma_{n-1}^2,
$$
$$
w_r = \sigma_{n-1} \ldots \sigma_{1} x_r \sigma_{1}^{-1} \ldots \sigma_{n-1}^{-1},~
r=1,\ldots,2g.
$$

\begin{prop} \label{theorem6}
The group $D_n(\Sigma_{g}),$ $n>1,$ admits a group presentation with  generators:
$$
\sigma_1, \ldots, \sigma_{n-2}, x_1, \ldots, x_{2g},
\tau_1, \ldots, \tau_{n-1},~w_1, \ldots, w_{2g},
$$
and  relations:

\noindent -- Braid relations:

\noindent (B1)~$\sigma_i \sigma_{i+1} \sigma_i = \sigma_{i+1} \sigma_i \sigma_{i+1},~~1 \leq i \leq n-3,$

\noindent (B2)~$\sigma_i \sigma_{j} = \sigma_{j} \sigma_i,~~|i - j| > 1,~~1 \leq i,j \leq n-2,$

\noindent (B3)~$\sigma_k^{-1} \tau_l \sigma_k = \tau_l,~k \neq l-1, l,$

\noindent (B4)~$\sigma_{l-1}^{-1} \tau_l \sigma_{l-1} = \tau_{l-1},$

\noindent (B5)~$\sigma_l^{-1} \tau_l \sigma_l = \tau_l \tau_{l+1} \tau_l^{-1},~l \neq n-1.$

\noindent -- Mixed relations:

\noindent (R1.1)~$x_r \sigma_i = \sigma_i x_r,$~~$2 \leq i \leq n-2,$~~$1 \leq r \leq 2g,$

\noindent (R1.2)~$x_r \tau_i = \tau_i x_r,$~~$2 \leq i \leq n-1,$

\noindent (R1.3)~$w_r \sigma_i = \sigma_i w_r,$~~$1 \leq r \leq 2g,$~~$1 \leq i \leq n-2,$

\noindent (R2.1)~$(\sigma_1^{-1} x_r \sigma_1^{-1}) x_r = x_r (\sigma_1^{-1} x_r \sigma_1^{-1}),$
~~$1 \leq r \leq 2g,$

\noindent (R2.2)~$x_r^{-1} w_r x_r = \tau_1^{-1} w_r \tau_1,$
~~$1 \leq r \leq 2g,$

\noindent (R2.3)~$x_r^{-1} \tau_1 x_r = \tau_1^{-1} w_r \tau_1 w_r^{-1} \tau_1,$
~~$1 \leq r \leq 2g,$

\noindent (R3.1)~$(\sigma_1^{-1} x_s \sigma_1) x_r = x_r (\sigma_1^{-1} x_s \sigma_1),$
~~$1 \leq s < r \leq 2g,$~~$(s, r) \neq (2m-1, 2m),$

\noindent (R3.2)~$x_r^{-1} (\tau_1^{-1} w_s \tau_1) x_r = \tau_1^{-1} w_s \tau_1,$
~~$1 \leq s < r \leq 2g,$~~$(s, r) \neq (2m-1, 2m),$

\noindent (R3.3)~$x_s w_r = w_r x_s,$
~~$1 \leq s < r \leq 2g,$~~$(s, r) \neq (2m-1, 2m),$

\noindent (R4.1)~$(\sigma_1^{-1} x_{2m-1} \sigma_1^{-1}) x_{2m} = x_{2m} (\sigma_1^{-1} x_{2m-1}
\sigma_1),$~~$1 \leq m \leq g,$

\noindent (R4.2)~$x_{2m}^{-1} ( \tau_1^{-1} w_{2m-1} ) x_{2m} = \tau_1^{-1} w_{2m-1} \tau_1,$
~~$1 \leq m \leq g,$

\noindent (R4.3)~$x_{2m-1}^{-1} w_{2m}  x_{2m-1} = \tau_1^{-1} w_{2m},$
~~$1 \leq m \leq g,$

\noindent (RT.1)~$[x_1^{-1}, x_2] [x_3^{-1}, x_4] \ldots [x_{2g-1}^{-1}, x_{2g}] =
\sigma_1 \sigma_2 \ldots \sigma_{n-3} \sigma_{n-2}^{2} \sigma_{n-3} \ldots \sigma_{2}
\sigma_1 \tau_1,$

\noindent (RT.2)~$[w_1^{-1}, w_2] [w_3^{-1}, w_4] \ldots [w_{2g-1}^{-1}, w_{2g}] =
\tau_1 \tau_2 \ldots \tau_{n-1}.$

\end{prop}

Consider the subgroup $V_{n-1,g}$ of $D_{n}(\Sigma_{g})$
generated by $\{ \tau_2, \ldots, \tau_{n-1}, w_1, $ $\ldots, w_{2g} \}$.
The group $V_{n-1,g}$  is free and normal in $D_{n}(\Sigma_{g})$~\cite{Bel}.
Using  relation (RT.2) one finds that
$$
\tau_1 = [w_1^{-1}, w_2] [w_3^{-1}, w_4] \ldots
[w_{2g-1}^{-1}, w_{2g}] \tau_{n-1}^{-1} \ldots \tau_{2}^{-1}
$$
and therefore $
\tau_1 \in V_n$.
 Also, consider the subgroup
$$
\overline{B}_{n-1}(\Sigma_{g}) = \langle \sigma_1, \ldots, \sigma_{n-2}, x_1,
\ldots, x_{2g} \rangle
$$
of $D_n(\Sigma_{g})$. The group $\overline{B}_{n-1}(\Sigma_{g})$ is not isomorphic to
$B_{n-1}(\Sigma_{g})$.

In fact, one can prove that the relations (B1), (B2), (R1.1), (R2.1), (R3.1), (R4.1) and (RT.1)
are defining relations of
$\overline{B}_{n-1}(\Sigma_g)$ and therefore we can deduce that 
$B_{n-1}(\Sigma_{g}) \simeq \overline{B}_{n-1}(\Sigma_{g})/ \langle\langle \tau_1 \rangle \rangle$. 
From the other relations we can see
that $\overline{B}_{n-1}(\Sigma_{g})$ acts on $V_n$ by
conjugacy. Hence, the following theorem holds.

\begin{prop} \label{theorem6b}
There exists a representation $\rho_V: \overline{B}_{n-1}(\Sigma_{g}) \to {\rm Aut} (V_{n-1,g})$, induced by the action
by conjugacy of  $\overline{B}_{n-1}(\Sigma_{g})$ on the free group $V_{n-1,g}=\langle \tau_2, \ldots, \tau_{n-1}, w_1, $ $\ldots, w_{2g} \rangle$ given algebraically as follows:
\\

\noindent --Action on generators $\tau_2, \ldots, \tau_{n-1}$:

\noindent (S1) $\tau_l^{\sigma_k} = \tau_l,$~$k \neq l-1, l,$
\\

\noindent (S2) $\tau_l^{\sigma_{l-1}} = \tau_{l-1},$
\\

\noindent (S3) $\tau_l^{\sigma_l} = \tau_{l+1}^{\tau_l^{-1}},$~$l \neq n-1,$
\\

\noindent (S4) $\tau_i^{x_r} = \tau_i,$~$2 \leq i \leq n-1,$~$1 \leq r \leq 2g,$
\\

\noindent --Action on generators $w_1, \ldots, w_{2g} $:

\noindent (S5) $w_r^{\sigma_i} = w_r,$~$1 \leq r \leq 2g,$~~$1 \leq  i \leq n-2,$
\\

\noindent (S6) $w_r^{x_r} = w_r^{\tau_1},$~$1 \leq r \leq 2g,$
\\

\noindent (S7) $\tau_1^{x_r} = \tau_1^{w_r^{-1} \tau_1},$~$1 \leq r \leq 2g,$
\\

\noindent (S8) $w_s^{x_r} = w_s^{[w_r^{-1}, \tau_1]},$~$1 \leq s < r \leq 2g,$~$(s, r) \neq (2m-1, 2m),$
\\

\noindent (S9) $w_r^{x_s} = w_r,$~$1 \leq s < r \leq 2g,$~$(s, r) \neq (2m-1, 2m),$
\\

\noindent (S10) $w_{2m-1}^{x_{2m}} = [\tau_1, w_{2m}^{-1}] w_{2m-1} \tau_1,$
\\

\noindent (S11) $w_{2m}^{x_{2m-1}} = \tau_1^{-1} w_{2m},$
\\
where
$$
\tau_1 = [w_1^{-1}, w_2] [w_3^{-1}, w_4] \ldots
[w_{2g-1}^{-1}, w_{2g}] \tau_{n-1}^{-1} \ldots \tau_{2}^{-1}.
$$
\end{prop}

From Proposition~\ref{theorem6b} we get a representation of $B_{n-1}(\Sigma_g)$ in
${\rm Aut} (V_{n-1,g})/\langle \langle t \rangle \rangle$, where
 $t \in {\rm Aut} (V_{n-1,g})$ is the conjugation by the word
$$
[w_1^{-1}, w_2] \ldots [w^{-1}_{2g-1}, w_{2g}]\tau_{n-1}^{-1} \ldots \tau_2^{-1}.
$$
Therefore, we have the following result:

\begin{prop} \label{theorem7}
There exists a representation $\rho_{\widetilde{V}}: B_{n-1}(\Sigma_{g}) \to {\rm Out} (V_{n-1,g})$,
given algebraically as follows:
\\

\noindent --Action (up to conjugacy) on generators $\tau_2, \ldots, \tau_{n-1}$:

\noindent (S1) $\tau_l^{\sigma_k} = \tau_l,$~$k \neq l-1, l,$
\\

\noindent (S2) $\tau_l^{\sigma_{l-1}} = \tau_{l-1},$
\\

\noindent (S3) $\tau_l^{\sigma_l} = \tau_{l+1}^{\tau_l^{-1}},$~$l \neq n-1,$
\\

\noindent (S4) $\tau_i^{x_r} = \tau_i,$~$2 \leq i \leq n-1,$~$1 \leq r \leq 2g,$
\\

\noindent --Action (up to conjugacy) on generators $w_1, \ldots, w_{2g} $:

\noindent (S5) $w_r^{\sigma_i} = w_r,$~$1 \leq r \leq 2g,$~~$1 \leq  i \leq n-2,$
\\

\noindent (S6) $w_r^{x_r} = w_r^{\tau_1},$~$1 \leq r \leq 2g,$
\\

\noindent (S7) $\tau_1^{x_r} = \tau_1^{w_r^{-1} \tau_1},$~$1 \leq r \leq 2g,$
\\

\noindent (S8) $w_s^{x_r} = w_s^{[w_r^{-1}, \tau_1]},$~$1 \leq s < r \leq 2g,$~$(s, r) \neq (2m-1, 2m),$
\\

\noindent (S9) $w_r^{x_s} = w_r,$~$1 \leq s < r \leq 2g,$~$(s, r) \neq (2m-1, 2m),$
\\

\noindent (S10) $w_{2m-1}^{x_{2m}} = [\tau_1, w_{2m}^{-1}] w_{2m-1} \tau_1,$
\\

\noindent (S11) $w_{2m}^{x_{2m-1}} = \tau_1^{-1} w_{2m},$
\\
where
$$
\tau_1 = [w_1^{-1}, w_2] [w_3^{-1}, w_4] \ldots
[w_{2g-1}^{-1}, w_{2g}] \tau_{n-1}^{-1} \ldots \tau_{2}^{-1}.
$$
\end{prop}

We don't know if the representation  $\rho_{\widetilde{V}}$ of $B_{n-1}(\Sigma_{g})$ is faithful or not.

\section{Artin-Tits groups}

We recall that classical braid groups are also called Artin-Tits groups of type $\mathcal{A}$. More precisely, let  $(W,S)$ be 
a Coxeter system and let us denote 
the order of the element $s t$ in $W$ by $m_{s, \, t}$ (for $s, t \in S$).
Let $A(W)$ be the group defined by the following group presentation:
$$
A(W)=\langle S \, | \underbrace{s t\cdots}_{m_{s, \, t}}= \underbrace{t s\cdots}_{m_{s, \, t}}\; \mbox{ for any $s\not=t \in S$
with $m_{s, \, t} < +\infty$} \, \rangle\;.
$$
The group $A(W)$ is the Artin-Tits group associated to $W$. The group
$A(W)$ is said to be of spherical type if  $W$ is finite. There exists three infinite families of finite Coxeter groups
usually denoted respectively as of type $\mathcal{A}$, $\mathcal{B}$ and $\mathcal{D}$.
Artin braid groups correspond to the family of Artin-Tits
groups associated to Coxeter groups of type $\mathcal{A}$ and braid groups of the annulus coincide with Artin-Tits
groups associated to Coxeter groups of type $\mathcal{B}$. In this Section we give a faithful representation 
of the remaining infinite family of spherical Artin-Tits groups, associated to Coxeter groups
of type $\mathcal{D}$.

First let us denote by $A(\mathcal{D}_n)$ the $n$-th   Artin-Tits group of type $\mathcal{D}$ with the group presentation
provided by its Coxeter graph (see Figure~1), where vertices $\delta_1, \dots, \delta_n$ are the generators of the group and two generators $\delta_i, \delta_j$ verify the relation $\delta_i \delta_j \delta_i= \delta_j \delta_i \delta_j$  if they are related by an edge and elsewhere they commute.

\begin{figure}[htp]
	\centering
	\includegraphics[width=5cm,height=2.3cm,bb=0 0 283 122]{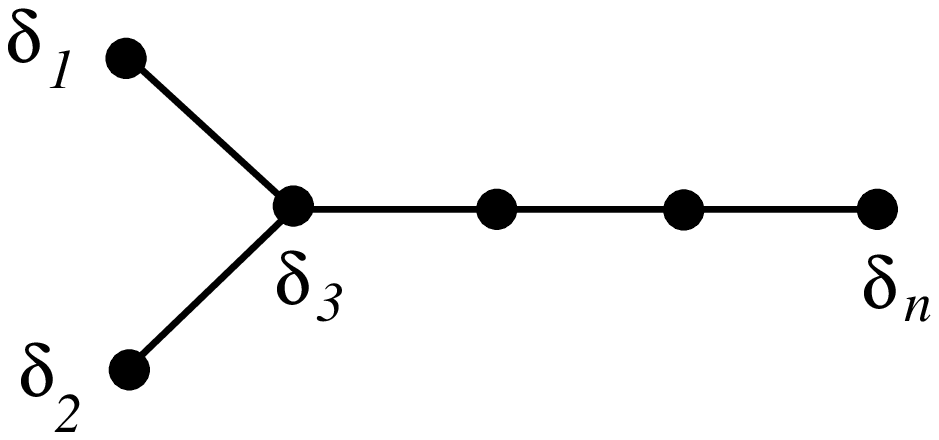}
	\label{fig:artinD}
\caption{}
\end{figure}

We denote by $\pi_D: A(\mathcal{D}_n) \to B_n$ the epimorphism defined by 
$\pi_D(\delta_1)=\pi_D(\delta_2)=\sigma_1$ and 
$\pi_D(\delta_i)=\sigma_{i-1}$ for $i=3, \ldots, n$, which admits the section
$s_D:    B_n \to  A(\mathcal{D}_n)$ defined by 
$s_D(\sigma_i)= \delta_{i+1}$ for $i=1, 2, \ldots, n-1$.

\begin{prop} \label{artind} \cite[Proposition 2.3]{CP}
\begin{enumerate}
\item The  representation $\rho_D: B_n \to Aut(F_{n-1})$ given algebraically by:

$$
\rho_D(\sigma_{1}) : 
\left\{ \begin{array}{ll}
x_{1} \to x_{1}, & \\
x_{j} \to x_{1}^{-1} x_j,  j \neq i,& 
\end{array} \right.
$$

and for $2 \le i \le n-1$,

$$ \rho_D(\sigma_{i}) : 
\left\{
\begin{array}{ll}
x_{i-1} \longmapsto x_{i}, & \\
x_{i} \longmapsto x_{i} x_{i-1}^{-1} x_i,  &\\
x_{j} \longmapsto x_{j}, &  j\neq i-1, i.
\end{array} \right.
$$
 
 is well defined and faithful. 
\item The group  $A(\mathcal{D}_n)$ is isomorphic to $F_{n-1} \rtimes_{\rho_D} B_n$, where the projection on the second factor is  $\pi_D$ and the section $B_n \to     F_{n-1} \rtimes_{\rho_D} B_n$     is just $s_D$. 
\item In particular, $\ker \pi_D=F_{n-1}$ is  freely generated by 
$\lambda_1, \ldots, \lambda_{n-1}$, where $\lambda_1= \delta_1 \delta_2^{-1}$ and $\lambda_i= (\delta_{i+1} \cdots \delta_3)(\delta_1 \delta_2^{-1})
(\delta_{i+1} \cdots \delta_3)^{-1}$ for
$i=2, \ldots, n-1$. 
\end{enumerate}
\end{prop}

The first item of Proposition~\ref{artind} was already established in  \cite{PV} by topological means and the third item is actually proven in the proof of Proposition 2.3 in \cite{CP}.

From Proposition~\ref{artind} we can deduce a faithful representation of $A(\mathcal{D}_n)$  into ${\rm Aut}(F_{n})$.

\begin{prop} \label{thmartind} 
The  representation $\iota: A(\mathcal{D}_n) \to {\rm Aut}(F_{n})$ given algebraically by:

$$
\iota(\delta_{1}) : \left\{
\begin{array}{ll}
x_{1} \longmapsto x_{1}, & \\
x_{j} \longmapsto x_{j} x_1^{-1},&   j \neq 1, \, n, \\
x_{n} \longmapsto x_1 x_{n} x_1^{-1},&\\
\end{array} \right.
$$

$$
\iota(\delta_{2}) : \left\{
\begin{array}{ll}
x_{1} \longmapsto x_{1}, & \\
x_{j} \longmapsto x_{1}^{-1} x_j, & j \neq 1, \, n, \\
x_{n} \longmapsto x_{n}, &\\
\end{array} \right.
$$

and for $3 \le i \le n$,

$$ \iota(\delta_{i}) : \left\{
\begin{array}{ll}
x_{i-2} \longmapsto x_{i-1}, & \\
x_{i-1} \longmapsto x_{i-1} x_{i-2}^{-1} x_{i-1},  &\\
x_{j} \longmapsto x_{j}, &  j\neq i-2, i-1.
\end{array} \right.
$$ 
 is well defined and faithful.
\end{prop}
\begin{proof}
From  Proposition~\ref{artind} it follows that   $F_{n-1} \rtimes_{\rho_D} B_n$
is isomorphic to $A(\mathcal{D}_n)$ through the morphism 
$\chi: F_{n-1} \rtimes_{\rho_D} B_n \to A(\mathcal{D}_n)$ defined algebraically as follows:
$\chi(x_i)=\lambda_i$
for any generator $x_i$ of $F_{n-1}$ and $\chi(\sigma_i)=\delta_{i+1}$ for any generator $\sigma_i$
of $B_{n}$. 

Remark that $F_{n-1} \rtimes_{\rho_D} B_n$ is a subgroup of $F_{n-1} \rtimes {\rm Aut}(F_{n-1})$.
Now, considering the natural inclusion of ${\rm Aut}(F_{n-1})$ into  ${\rm Aut}(F_{n})$ leaving the generator
$x_n$ invariant and the action by conjugacy of $F_{n-1}$ on $F_n$, we obtain a morphism
$\phi: F_{n-1} \rtimes_{\rho_D} B_n \to {\rm Aut}(F_{n})$, which can be easily proved to be injective. 
  
By direct 
calculation one can verify that $\iota(\lambda_{i})$ is the conjugacy by $x_i$ for $i=1, \ldots, n-1$;
moreover, since $\iota \circ \chi$ restricted to $B_n$ coincides with the representation of
 $B_n$ into ${\rm Aut}(F_{n})$ obtained composing  the map $\rho_D$ with 
 the natural inclusion of ${\rm Aut}(F_{n-1})$ into  ${\rm Aut}(F_{n})$ leaving the generator
$x_n$ invariant, we can deduce that
 $\iota \circ \chi= \phi$ and therefore  the  
claim follows.
\end{proof}


\end{document}